\documentclass[11pt]{amsart}
\usepackage{CJK}
\usepackage{amssymb}
\usepackage{amsmath}
\usepackage{amsthm}
\usepackage{bm}

\usepackage{graphicx}

\newtheorem{Prop}{Proposition}[section]
\newtheorem{Thm}[Prop]{Theorem}
\newtheorem{Lem}[Prop]{Lemma}
\newtheorem{Cor}[Prop]{Corollary}
\newtheorem{Conj}[Prop]{Conjecture}

\newtheorem*{Thm0}{Theorem}

\theoremstyle{definition}
\newtheorem{Ex}[Prop]{Example}
\newtheorem{Def}[Prop]{Definition}

\begin{document}

\title{The genus filtration in the smooth concordance group}

\author{Shida Wang}

\address{Department of Mathematics, Indiana University, Bloomington, IN 47405}
\email{wang217@indiana.edu}

\begin{abstract}

We define a filtration of the smooth concordance group based on the genus of representative knots.
We use the Heegaard Floer epsilon and Upsilon invariants to prove the quotient groups with respect to this filtration are infinitely generated.
Results are applied to three infinite families of topologically slice knots.

\end{abstract}

\maketitle

\section{Introduction}

Let $\mathcal{C}$ be the smooth concordance group.
Let $\mathcal{G}_k$ denote the subgroup of $\mathcal{C}$ generated by knots of genus not greater than $k$.
Clearly $\mathcal{G}_0\subset\mathcal{G}_1\subset\mathcal{G}_2\subset\cdots\subset\mathcal{G}_k\subset\cdots\subset\mathcal{C}$ and $\bigcup_{k=1}^\infty=\mathcal{C}$.
This gives a filtration in $\mathcal{C}$. We call it the \emph{genus filtration}.

There is another way to understand $\mathcal{G}_k$. Recall that the concordance genus $g_c$ of a knot $K$ is defined to be the minimal genus of a knot $K'$ concordant to $K$.
It is obvious that $g_c(K)\\=\min\{k|K\text{ is concordant to }K'\text{ and }g(K')\leqslant k\}$.
This motivates the following definition.

\begin{Def}The \emph{splitting concordance genus} of a knot $K$ is $$g_{\mathrm{sp}}(K):=\min\{k|K\text{ is concordant to }K_1\#\cdots\#K_m\text{ and }g(K_1),\cdots,g(K_m)\leqslant k\}.$$\end{Def}

That is to say, $g_{\mathrm{sp}}(K)$ is the filtration level of $K$ in $\mathcal{G}_0\subset\mathcal{G}_1\subset\mathcal{G}_2\subset\cdots\subset\mathcal{G}_k\subset\cdots\subset\mathcal{C}$.

By \cite{pretzel}, $\mathcal{G}_1$ contains a $\mathbb{Z}^\infty$ subgroup whose elements are topologically slice.

Let $\mathcal{C}_{TS}\subset\mathcal{C}$ be the subgroup of topologically slice knots.
Recently several results revealing that the group $\mathcal{C}_{TS}$ is quite large have been proven.
For example, \cite{1nn1} and \cite{upsilon} show $\mathcal{C}_{TS}$ contains a $\mathbb{Z}^\infty$ direct summand;
\cite{Alexander1} shows $\mathcal{C}_{TS}$ contains a $\mathbb{Z}^\infty$ subgroup whose nonzero elements are not concordant to Alexander polynomial one knots;
\cite{order2} shows $\mathcal{C}_{TS}$ contains a $\mathbb{Z}_2^\infty$ subgroup whose nonzero elements are not concordant to Alexander polynomial one knots.

We will show $\mathcal{C}_{TS}$ is large in another sense.
We will prove $\mathcal{C}_{TS}\not\subset\mathcal{G}_k$ for any $k$.
Moreover, the difference between $\mathcal{C}_{TS}$ and $\mathcal{G}_k$ is large.
Corollary \ref{main} states that $\mathcal{C}_{TS}/(\mathcal{C}_{TS}\cap\mathcal{G}_{k})$ contains a direct summand isomorphic to $\mathbb{Z}^\infty$.

Our examples will be built from those of \cite{epsilon} and \cite{upsilon}.
Let $Wh(K)$ denote the untwisted Whitehead double of a knot $K$ and $K_{p,q}$ denote the $(p,q)$-cable of $K$.
Let $J_n=(Wh(T_{2,3}))_{n,n+1}\#-T_{n,n+1}$ and let $J_n'=(Wh(T_{2,3}))_{n,2n-1}\#-T_{n,2n-1}$.
These knots are topologically slice and used to prove the following theorems.

\begin{Thm0} \emph{(\cite[Theorem 1]{1nn1})} The group $\mathcal{C}_{TS}$ contains a summand isomorphic to $\mathbb{Z}^\infty$ generated by $\{J_n\}_{n=2}^\infty$.\end{Thm0}

\begin{Thm0} \emph{(\cite[Theorem 1.20]{upsilon})} The topologically slice knots $\{J_n'\}_{n=2}^\infty$ form a basis for a free direct summand of $\mathcal{C}_{TS}$.\end{Thm0}

We will prove the following results.

\begin{Thm}\label{Hom} $\{J_n\}_{n=k}^\infty$ form a basis for a $\mathbb{Z}^\infty$ summand of $\mathcal{C}_{TS}/(\mathcal{C}_{TS}\cap\mathcal{G}_{\lfloor\frac{k}{2}\rfloor})$ for any $k\geqslant2$.\end{Thm}

\begin{Thm}\label{OSS} $\{J_n'\}_{n=k}^\infty$ form a basis for a $\mathbb{Z}^\infty$ summand of $\mathcal{C}_{TS}/(\mathcal{C}_{TS}\cap\mathcal{G}_{k-1})$ for any $k\geqslant2$.\end{Thm}

Hence we have

\begin{Cor}\label{main} For any $k\in\mathbb{N}$, $\mathcal{C}_{TS}\not\subset\mathcal{G}_k$.
Moreover, the quotient group $\mathcal{C}_{TS}/(\mathcal{C}_{TS}\cap\mathcal{G}_{k})$ contains a direct summand isomorphic to $\mathbb{Z}^\infty$.\end{Cor}

One can define another subgroup $\mathcal{H}_k$ of $\mathcal{C}$ to be generated by knots of 4--genus not greater than $k$.
Clearly $\mathcal{G}_k\subset\mathcal{H}_k$.
It is natural to ask whether $\mathcal{H}_k/\mathcal{G}_k$ is infinitely generated.
We show the answer is affirmative by proving

\begin{Thm}\label{g4} The quotient group $\mathcal{C}_{TS}/(\mathcal{C}_{TS}\cap\mathcal{G}_{k})$ contains a subgroup
isomorphic to $\mathbb{Z}^\infty$ whose basis elements have slice genus 1 for any $k\geqslant2$.\end{Thm}

\begin{Conj}\ \\(1) For any $k\in\mathbb{N}$,
$(\mathcal{C}_{TS}\cap\mathcal{G}_{k+1})/(\mathcal{C}_{TS}\cap\mathcal{G}_{k})$ contains a direct summand isomorphic to $\mathbb{Z}^\infty$ whose basis elements have slice genus 1.

(2) For any $k\in\mathbb{N}$,
$\mathcal{C}/\mathcal{H}_n$ is nontrivial.\end{Conj}

This paper is organized as follows. Section 2 uses Alexander polynomials to prove the splitting concordance genus can be arbitrarily large.
Section 3 reviews Hom's $\varepsilon$ invariant and develops an obstruction, which is used to prove Theorems \ref{Hom} and \ref{g4}.
Section 4 develops an obstruction by the Upsilon invariant and proves Theorem \ref{OSS}.

\emph{Acknowledgments.} The author wishes to express sincere thanks to Professor Charles Livingston for proposing this study and carefully reading a draft of this paper.
Thanks also to Jennifer Hom for valuable comments.

\section{A first glance at the genus filtration}

The Alexander polynomial $\Delta_K(t)$ of a slice knot $K$ satisfies the Fox-Milnor condition that it factors as $t^{\pm n}f(t)f(t^{-1})$.
Also, the maximal exponent minus the minimal exponent, $\mathrm{breadth}\ \Delta_K(t)$, is a lower bound of twice genus for any knot $K$.
Based on these facts, we can prove the following theorem, generalizing \cite[Theorem 2.2]{concordance genus}.

\begin{Prop}For any knot $K$, suppose that $p(t)$ appears an odd number of times in the irreducible factorization of $\Delta_K(t)$ in $\mathbb{Z}[t,t^{-1}]$.
Then $g_{\mathrm{sp}}(K)\geqslant\frac{1}{2}\mathrm{breadth}\ p(t)$.\end{Prop}

\textbf{Proof.} By definition, $K$ is concordant to $K_1\#\cdots\#K_m$ and $g(K_1),\cdots,g(K_m)\leqslant g_{\mathrm{sp}}(K)$.
Thus $K\#-K_1\#\cdots\#-K_m$ is a slice knot
and its Alexander polynomial $\Delta_K(t)\Delta_{K_1}(t)\cdots\Delta_{K_m}(t)$ must factor as $t^{\pm n}f(t)f(t^{-1})$ for some $f\in\mathbb{Z}[t,t^{-1}]$.
If some $p(t)$ appear odd times in the irreducible factorization of $\Delta_K(t)$, it must appear in the irreducible factorization of one of $\Delta_{K_1}(t),\cdots,\Delta_{K_m}(t)$.
Since $\mathrm{breadth}\ \Delta_{K_1}(t),\cdots,\mathrm{breadth}\ \Delta_{K_m}(t)\leqslant2g_{\mathrm{sp}}(K)$, we get $2g_{\mathrm{sp}}(K)\geqslant\mathrm{breadth}\ p(t)$.\hfill$\Box$

\begin{Ex}The Alexander polynomial of the torus knot $T_{p,q}$ is $\Delta_{T_{p,q}}(t)=\frac{(t^{pq}-1)(t-1)}{(t^p-1)(t^q-1)}$,
in whose irreducible factorization the cyclotomic polynomial $\Phi_{pq}$ appears exactly once.
Hence $g_{\mathrm{sp}}(T_{p,q})\geqslant\frac{\varphi(pq)}{2}$, where $\varphi$ is Euler's totient function.
If $p$ and $q$ are prime, we have $g_{\mathrm{sp}}(T_{p,q})\geqslant\frac{(p-1)(q-1)}{2}$ and this is actually an equality because $g(T_{p,q})=\frac{(p-1)(q-1)}{2}$.\end{Ex}

\begin{Cor}$\mathcal{C}/\mathcal{G}_{k}$ is nontrivial for any $k\in\mathbb{N}$.\end{Cor}

Working a little harder, we can show the following.

\begin{Prop}$\mathcal{C}/\mathcal{G}_{k}$ contains an infinitely generated free subgroup for any $k\in\mathbb{N}$.\end{Prop}

\textbf{Proof.} Let $\{p_n\}_{n=1}^\infty$ be a sequence of strictly increasing prime numbers with $p_1>k$.
We will prove that the torus knots $\{T_{p_{2n-1},p_{2n}}\}_{n=1}^\infty$ are linearly independent in $\mathcal{C}/\mathcal{G}_{k}$.

Suppose towards a contradiction that $\#_{i=1}^l c_iT_{p_{_{2n_i-1}},p_{_{2n_i}}}$, where $0<n_1<\cdots<n_l$ and the $c_i$'s are nonzero integers,
is concordant to $K_1\#\cdots\#K_m$ with $g(K_1),\cdots,g(K_m)\leqslant k$.
Notice that $\Delta_{T_{p_{_{2n_i-1}},p_{_{2n_i}}}}(t)=\Phi_{p_{_{2n_i-1}}p_{_{2n_i}}}$,
where $\Phi_{p_{_{2n_i-1}}p_{_{2n_i}}}$ is the cyclotomic polynomial,
which is irreducible of degree $(p_{_{2n_i-1}}-1)(p_{_{2n_i}}-1)$.
By a combinatorial formula \cite[Proposition 1]{signature} for the Tristram-Levine signature functions of torus knots,
$\sigma_{\omega}(T_{p_{_{2n_i-1}},p_{_{2n_i}}})$ jumps by $\pm2$ at the primitive $p_{_{2n_i-1}}p_{_{2n_i}}$th roots of unity.
Since the $p_{_{2n_i-1}}p_{_{2n_i}}$'s are distinct for $i=1,\cdots,l$,
we know $\sigma_{\omega}(\#_{i=1}^l c_iT_{p_{_{2n_i-1}},p_{_{2n_i}}})$ has a jump discontinuity at a primitive $p_{_{2n_1-1}}p_{_{2n_1}}$th root of unity.
Hence $\sigma_{\omega}(K_1\#\cdots\#K_m)$ also has a jump discontinuity at a primitive $p_{_{2n_1-1}}p_{_{2n_1}}$th root of unity,
and so does one of $\sigma_{\omega}(K_1),\cdots,\sigma_{\omega}(K_m)$.
Without loss of generality, assume that $\sigma_{\omega}(K_1)$ has a jump discontinuity at a primitive $p_{_{2n_1-1}}p_{_{2n_1}}$th root of unity.
Since jump discontinuities of the Tristram-Levine signature function can only appear at roots of the Alexander polynomial,
it follows that $\Delta_{K_1}(t)$ has a root at a primitive $p_{_{2n_1-1}}p_{_{2n_1}}$th root of unity and thus is divisible by $\Phi_{p_{_{2n_i-1}}p_{_{2n_i}}}$,
but this is impossible because $\mathrm{deg}\Delta_{K_1}(t)\leqslant2g(K_1)\leqslant2k<(p_{_{2n_i-1}}-1)(p_{_{2n_i}}-1)$.
\hfill$\Box$

\section{Obstruction by $\varepsilon$ invariant}

We assume the reader is familiar with knot Floer homology, defined by Ozsv\'{a}th-Szab\'{o} \cite{CFKdef1} and independently Rasmussen \cite{CFKdef2},
and the $\varepsilon$ invariant, defined by Hom \cite{epsilon}.
We briefly recall some of their properties for later use.

\subsection{The knot Floer complex and $\varepsilon$ invariant}

To a knot $K\subset S^3$, a doubly filtered, free, finitely generated chain complex over $\mathbb{F}[U,U^{-1}]$, denoted by $CFK^\infty(K)$, is associated,
where $\mathbb{F}$ is the field with two elements.
The two filtrations are called algebraic and Alexander filtrations and the grading of the chain complex is called the homological or Maslov grading.
Multiplication by $U$ shifts each filtration down by one and lowers the homological grading by two.
$CFK^\infty(K)$ is an invariant of $K$ up to filtered chain homotopy equivalence.
Furthermore, up to filtered chain homotopy equivalence, one can assume the differential strictly lowers at least one of the filtrations \cite{CFKdef2}.

A quick corollary from \cite[Theorem 1.2]{CFKrange} is the following.

\begin{Prop}\label{range} If $K$ has genus g, then all elements of $CFK^\infty(K)$ have filtration level $(i,j)$ where $-g\leqslant i-j\leqslant g$.\end{Prop}

For appropriate $S\subset\mathbb{Z}\oplus\mathbb{Z}$, let $C\{S\}$ denote the corresponding subquotient complex of $CFK^\infty(K)$.
For example, $C\{i=0\}=CFK^\infty(K)\{i\leqslant0\}/CFK^\infty(K)\{i<0\}$.
The invariant $$\tau(K)=\min\{s|C\{i=0,j\leqslant s\}\rightarrow C\{i=0\}\text{ induces an nontrivial map on homology}\}$$ is proven to be a smooth concordance invariant in \cite{tau}.

For any knot $K$, Hom defines an invariant called $\varepsilon$ taking on values $-1$, 0 or 1 \cite{epsilon}, which has the following properties.

\begin{Prop}\label{epsilon properties} \emph{(\cite[Proposition 3.6]{epsilon})} $\varepsilon$ satisfies the following properties:

(1) if $K$ is smoothly slice, then $\varepsilon(K)=0$;

(2) $\varepsilon(-K)=-\varepsilon(K)$;

(3) if $\varepsilon(K)=\varepsilon(K')$, then $\varepsilon(K\#K')=\varepsilon(K)=\varepsilon(K')$;

(4) if $\varepsilon(K)=0$, then $\varepsilon(K\#K')=\varepsilon(K')$.\end{Prop}

Thus the relation $\sim$, defined by $K\sim K'\Leftrightarrow\varepsilon(K\#-K')=0$, is an equivalence relation coarser than the smooth concordance.
It gives an equivalence relation on $\mathcal{C}$ called $\varepsilon$-equivalence.
The $\varepsilon$-equivalence class of $K$ is denoted by $[\![K]\!]$.
All $\varepsilon$-equivalence classes form a group $\mathcal{F}$ (also denoted by $\mathcal{CFK}$ in \cite{1nn1}), which is a quotient group of $\mathcal{C}$.
The kernel of the natural homomorphism from $\mathcal{C}$ to $\mathcal{F}$ is $\{[K]\in\mathcal{C}|\varepsilon(K)=0\}$.

According to \cite[Proposition 4.1]{ordered}, $\varepsilon$ induces a total order on $\mathcal{F}$.
The proof uses Proposition \ref{epsilon properties}.
The total order is defined by $[\![K]\!]>[\![K']\!]\Leftrightarrow\varepsilon(K\#-K')=1$.
Moreover, this order respects the addition operation on $\mathcal{F}$.
Therefore there is a quotient homomorphism from $\mathcal{C}$ to the totally ordered abelian group $\mathcal{F}$, which can be used to show linear independence in $\mathcal{C}$.

\subsection{Some facts about totally ordered abelian group}

Let $G$ be a totally ordered abelian group, that is, an abelian group with a total order respecting the addition operation. Denote its identity element by 0.

The absolute value of an element $a\in G$ is defined to be $|a|=\left\{\begin{aligned}&a&\text{ if }a\geqslant0,\\&-a&\text{ if }a<0.\end{aligned}\right.$

\begin{Def} Two nonzero elements $a$ and $b$ of $G$ are \emph{Archimedean equivalent}, denoted by $a\sim_A b$,
if there exists a natural number $N$ such that $N\cdot|a|>|b|$ and $N\cdot|b|>|a|$.
If $a$ and $b$ are not Archimedean equivalent and $|a|<|b|$, we say that $b$ dominates $a$ and we write $|a|\ll|b|$.\end{Def}

\begin{Def} (Property A) An element $a\in G$ \emph{satisfies Property A} if for every $b\in G$ such that $b\sim_A a$, we have that $b=ka+c$,
where $k$ is an integer and $c$ is dominated by $a$.\end{Def}

We have the following two facts. See \cite[Lemma 4.7]{ordered} and \cite[Proposition 1.3]{1nn1} for their proofs, respectively.

\begin{Lem}\label{domination implies independence} If $0<a_1\ll a_2\ll a_3\ll\cdots$ in $G$, then $a_1,a_2,a_3,\cdots$ are linearly independent in $G$.\end{Lem}

\begin{Lem}\label{PropertyA implies summand} If $0<a_1\ll a_2\ll a_3\ll\cdots$ in $G$ and each $a_i$ satisfies Property A,
then $a_1,a_2,a_3,\cdots$ generate (as a basis) a direct summand isomorphic to $\mathbb{Z}^\infty$ in $G$.\end{Lem}

\cite{ordered} shows

\begin{Lem}\label{1nn1domination} $0<[\![J_n]\!]\ll[\![J_{n+1}]\!]$ for any $n\geqslant2$.\end{Lem}

\cite{1nn1} shows

\begin{Lem}\label{1nn1PropertyA} $[\![J_n]\!]$ satisfies Properties A for any $n\geqslant2$.\end{Lem}

It is straightforward to check that $\{a\big||a|\ll x\}$ is a subgroup of $G$ for any $x>0$ in $G$. Denote this subgroup by $G_x$.
Let $\varphi_x$ be the quotient homomorphism.

\begin{Prop}\label{quotient} $G/G_x$ is a totally ordered abelian group with $\varphi_x(a)<\varphi_x(b)$ in $G/G_x$ defined by $a<b$ and $b-a\not\in G_x$.
If $0<a\ll b$ in $G$ and $b\not\in G_x$, then $0\leqslant\varphi_x(a)\ll\varphi_x(b)$ in $G/G_x$.
If $a$ satisfies Property A in $G$, then $\varphi_x(a)$ satisfies Property A in $G/G_x$.\end{Prop}

\textbf{Proof.} First we check that the relation $<$ in $G/G_x$ is well defined.
Suppose $\varphi_x(a)<\varphi_x(b)$. Let $c\in G_x$. We must show $\varphi_x(a+c)<\varphi_x(b)$ and $\varphi_x(a)<\varphi_x(b+c)$.
Since $b-a>0$ and $b-a\not\in G_x$, it is easy to verify that $b-a\gg|y|$ for any $y\in G_x$.
Thus $b-a\pm c>0$. Also $b-a\not\in G_x$ implies $b-a\pm c\not\in G_x$.
Hence $\varphi_x(a+c)<\varphi_x(b)$ and $\varphi_x(a)<\varphi_x(b+c)$,
which means the definition does not depend on the choices of $a$ and $b$.

Next we verify ``$<$" is a strict total order on $G/G_x$ that respects the addition operation.
For trichotomy, let $\varphi_x(a)$ and $\varphi_x(b)$ be two distinct elements in $G/G_x$.
Then $b-a\not\in G_x$. Thus $b-a\neq0$ and exactly one of $a<b$ and $b<a$ is true.
Hence exactly one of $\varphi_x(a)<\varphi_x(b)$ and $\varphi_x(b)<\varphi_x(a)$ is true by definition.
For transitivity, let $\varphi_x(a),\varphi_x(b),\varphi_x(c)\in G/G_x$ satisfy $\varphi_x(a)<\varphi_x(b)$ and $\varphi_x(b)<\varphi_x(c)$.
Then $a<b,b<c$ and $b-a,c-b\not\in G_x$. Immediately $a<c$.
Suppose towards a contradiction that $c-a\in G_x$. Then the fact that $b-a\gg|y|$ for any $y\in G_x$ implies $b-a-(c-a)>0$, which contradicts $b<c$.
Hence $c-a\not\in G_x$ and $\varphi_x(a)<\varphi_x(c)$ by definition.
For consistency with the addition operation, let $\varphi_x(a),\varphi_x(b),\varphi_x(c)\in G/G_x$ and $\varphi_x(a)<\varphi_x(b)$.
Then $a<b$ and $b-a\not\in G_x$. Thus $a+c<b+c$ and $(b+c)-(a+c)\not\in G_x$.
Hence $\varphi_x(a)+\varphi_x(c)=\varphi_x(a+c)<\varphi_x(b+c)=\varphi_x(b)+\varphi_x(c)$ by definition.

Finally we prove the domination and Property A in $G$ implies that in $G/G_x$.
Suppose $0<a\ll b$ in $G$ and $b\not\in G_x$. Then $0<Na<b$ for any $N\in\mathbb{N}$.
Also, the fact that $b\gg|y|$ for any $y\in G_x$ implies $Na+y<b,\forall y\in G_x$.
It follows that $b-Na>0$ and that $b-Na\not\in G_x$.
Hence $0\leqslant\varphi_x(a)\ll\varphi_x(b)$ in $G/G_x$ by definition.
Suppose $a$ satisfies Property A in $G$, that is, if $b\sim_A a$ in $G$ then $b=ka+c$ for some integer $k$ and some $c\in G$ dominated by $a$.
Without loss of generality we assume $\varphi_x(a)\neq0$.
Let $\varphi_x(b)\sim_A\varphi_x(a)$. So $b\sim_A a$ in $G$.
Otherwise either $|a|\ll|b|$ or $|b|\ll|a|$, which would imply $|\varphi_x(a)|\ll|\varphi_x(b)|$ or $|\varphi_x(b)|\ll|\varphi_x(a)|$
Thus $b=ka+c$ for some integer $k$ and some $c\in G$ dominated by $a$.
Thus $\varphi_x(b)=k\varphi_x(a)+\varphi_x(c)$.
Since $c$ is dominated by $a$, we know $\varphi_x(c)$ is dominated by $\varphi_x(a)$.
Hence $\varphi_x(a)$ satisfies Property A in $G/G_x$.
\hfill$\Box$

\subsection{Restriction on the Archimedean equivalence class by genus}

For any knot $K$ with $\varepsilon(K)=1$, a tuple of numerical invariant $\bm{a}^+(K)=(a_1(K),\cdots,a_n(K))$,
where $n$ depends on $K$ and $a_1(K),\cdots,a_n(K)$ are positive integers, is defined in \cite[Section 3]{1nn1}.
It is shown that $\bm{a}^+(K)$ is an invariant of the $\varepsilon$-equivalence class $[\![K]\!]$ (see \cite[Proposition 3.1]{1nn1}).

Computations in \cite{ordered} show

\begin{Lem}\label{1nn1a1a2} $\bm{a}^+(J_p)=(1,p,\cdots)$.\end{Lem}

$a_1$ and $a_2$ are useful in determining domination.

\begin{Lem}\label{a1a2greater} \emph{(\cite[Lemma 6.3 and 6.4]{ordered})} If $\bm{a}^+(K)=(a_1(K),\cdots)$ and $\bm{a}^+(K')=(a_1(K'),\cdots)$ with $a_1(K)>a_1(K')>0$,
then $[\![K]\!]\ll[\![K']\!]$.

If $\bm{a}^+(K)=(a_1(K),a_2(K),\cdots)$ and $\bm{a}^+(K')=(a_1(K'),a_2(K'),\cdots)$ with $a_1(K)=a_1(K')>0$ and $a_2(K)>a_2(K')>0$,
then $[\![K]\!]\gg[\![K']\!]$.\end{Lem}

Based on Proposition \ref{range}, the following is shown.

\begin{Lem}\label{a1a2tau} \emph{(\cite[Theorem 1.2 and Lemma 2.3]{gamma})} Suppose that $\varepsilon(K)=1$, and $a_2(K)$ is defined.
Then $|\tau(K)-a_1(K)-a_2(K)|\leqslant g(K)$.\end{Lem}

Next we prove our obstruction theorem.

\begin{Prop}\label{epsilon obstruction} Suppose $J$ is a knot with $\bm{a}^+(J)=(1,b,\cdots)$ with $b\geqslant2n$ for some positive integer $n$.
Then for any knot $K\in\mathcal{G}_n$, we have $|[\![K]\!]|\ll[\![J]\!]$.\end{Prop}

\textbf{Proof.} Before proving the proposition for $K\in\mathcal{G}_n$, first consider the case $g(K)\leqslant n$.
We may further assume that $[\![K]\!]>0$, since $[\![-K]\!]>0$ if $[\![K]\!]<0$ and it is trivial if $[\![K]\!]=0$.
Notice that $a_1(K)$ is always defined.
If $a_1(K)>1$, then $[\![K]\!]\ll[\![J]\!]$ by Lemma \ref{a1a2greater}.
If $a_1(K)=1$, then by \cite[Lemma 3.7]{1nn1} $a_2(K)$ is defined.
Observe that $\tau(K)-a_1(K)-a_2(K)\geqslant-g(K)$ by Lemma \ref{a1a2tau}.
Combining with $\tau(K)\leqslant g_4(K)\leqslant g(K)$, one has $g(K)-a_1(K)-a_2(K)\geqslant-g(K)$.
This implies $a_2(K)\leqslant2n-1$ if $a_1(K)=1$.
Hence $|[\![K]\!]|\ll[\![J]\!]$ by Lemma \ref{a1a2greater}.

Generally, let $K\in\mathcal{G}_n$. Then $K=K_1+\cdots+K_m$, where $g(K_i)\leqslant n$ for $i=1,\cdots,m$.
Since $[\![K]\!]=[\![K_1]\!]+\cdots+[\![K_m]\!]$, we know $|[\![K]\!]|\leqslant|[\![K_1]\!]|+\cdots+|[\![K_m]\!]|$.
Then the conclusion follows from the last paragraph.\hfill$\Box$

\subsection{Applying the obstruction to concrete families of knots}

\

\textbf{Proof of Theorem \ref{Hom}.} Fix an integer $k\geqslant2$.
Under the quotient homomorphism from $\mathcal{C}$ to $\mathcal{F}$,
the image of $\mathcal{G}_{\lfloor\frac{k}{2}\rfloor}$ is included in $\mathcal{F}_{[\![J_k]\!]}=\{[\![K]\!]\big||[\![K]\!]|\ll[\![J_k]\!]\}$
by Proposition \ref{epsilon obstruction} and Lemma \ref{1nn1a1a2}.
This gives a homomorphism from $\mathcal{C}/\mathcal{G}_{\lfloor\frac{k}{2}\rfloor}$ to $\mathcal{F}/\mathcal{F}_{[\![J_k]\!]}$.
By Lemma \ref{1nn1domination} and Proposition \ref{quotient}, the family $\{J_n\}_{n=k}^\infty$ maps to a family of elements with Property A and each term is dominated by the next.
Hence $\{J_n\}_{n=k}^\infty$ forms a basis of a direct summand isomorphic to $\mathbb{Z}^\infty$ by Lemma \ref{PropertyA implies summand}.
Note that since $J_n$'s are topologically slice, the above argument can be restricted to the subgroup
$\mathcal{C}_{TS}/(\mathcal{C}_{TS}\cap\mathcal{G}_{\lfloor\frac{k}{2}\rfloor})$ of $\mathcal{C}/\mathcal{G}_{\lfloor\frac{k}{2}\rfloor}$ to complete the proof.
\hfill$\Box$

\textbf{Proof of Theorem \ref{g4}.} Instead of $\{J_n\}$, we use another family of topologically slice knots $\{L_n\}$, where $L_n=(Wh(T_{2,3}))_{n,1}\#-(Wh(T_{2,3}))_{n-1,1}$.
These knots have slice genus 1 \cite[Lemma 3.1]{gamma}.
Also, \cite{gamma} computes that $a_1(L_n)=1$ and $a_2(L_n)=n$.
By the same argument as the above proof, except for applying Lemma \ref{domination implies independence} rather than Lemma \ref{PropertyA implies summand},
we immediately know $\{L_n\}_{n=2k}^\infty$ are linearly independent in $\mathcal{C}_{TS}/(\mathcal{C}_{TS}\cap\mathcal{G}_k)$.
\hfill$\Box$

\section{Obstruction by $\Upsilon$ invariant}

We refer to \cite{upsilon} for the definition and basic properties of $\Upsilon$ invariant.

For any knot $K$, $\Upsilon_K(t)$ is a piecewise linear function on $[0,2]$ whose derivative has finitely many discontinuities (\cite[Proposition 1.4]{upsilon}).
Thus, one can define $\Delta\Upsilon'_K(t_0)=\lim\limits_{t\rightarrow t_0^+}\Upsilon'_K(t)-\lim\limits_{t\rightarrow t_0^-}\Upsilon'_K(t)$ for any $t_0\in(0,2)$.

As an example, \cite[Proposition 1.4]{upsilon} computes the family $\{J_n'\}$ and get
$$\Delta\Upsilon'_{J_n'}(t)=\left\{\begin{aligned}&0&\text{ for }t<\frac{2}{2n-1},\\&2n-1&\text{ for }t=\frac{2}{2n-1}.\end{aligned}\right.$$

\cite[Corollary 1.12]{upsilon} shows $\Upsilon$ gives a homomorphism from $\mathcal{C}$ to the vector space of continuous functions on $[0,2]$.
Also, $K\mapsto\left\{\begin{aligned}\frac{1}{q}\Delta\Upsilon'_K(\frac{p}{q})\text{ if }p\text{ is even}\\
\frac{1}{2q}\Delta\Upsilon'_K(\frac{p}{q})\text{ if }p\text{ is odd}\end{aligned}\right.$
gives a homomorphism from $\mathcal{C}$ to $\mathbb{Z}$ for any $\frac{p}{q}\in(0,2)\cap\mathbb{Q}$.

The location of singularities of $\Upsilon$ is related to the genus of the knot, as in the following proposition.
Just like Lemma \ref{a1a2tau}, the proof is based on the fact in Proposition \ref{range} as well.

\begin{Prop}\label{denominator obstructs genus} \emph{(\cite[Theorem 8.2]{denominator})} Suppose that $\Delta\Upsilon'_K(t)$ is nonzero at $t=\frac{p}{q}$ with $\gcd(p,q)=1$.
Then $q\leqslant g(K)$ if $p$ is odd, and $q\leqslant2g(K)$ if $p$ is even.\end{Prop}

With this proposition, we can easily prove our obstruction theorem.

\begin{Prop}\label{upsilon obstruction} Suppose $K\in\mathcal{G}_n$ for some positive integer $n$. Then $\Delta\Upsilon'_K(t)=0$ for $t\in(0,\frac{1}{n})\cap\mathbb{Q}$.\end{Prop}

\textbf{Proof.} Before proving the proposition for $K\in\mathcal{G}_n$, first consider the case $g(K)\leqslant n$.
If $\Upsilon_K(t)$ has a singularity at a rational number $\frac{p}{q}$ with $\gcd(p,q)=1$, then Proposition \ref{denominator obstructs genus} implies $\frac{p}{q}\geqslant\frac{1}{n}$.

Generally, let $K\in\mathcal{G}_n$. Then $K=K_1+\cdots+K_m$, where $g(K_i)\leqslant n$ for $i=1,\cdots,m$.
If $\Upsilon_K(t)$ has a singularity at a rational number $\frac{p}{q}$, then so does one of $\Upsilon_{K_1}(t),\cdots,\Upsilon_{K_m}(t)$, since $\Upsilon$ is a homomorphism.
The conclusion follows from the last paragraph.\hfill$\Box$

\textbf{Proof of Theorem \ref{OSS}.} Fix an integer $k\geqslant2$.
If $K\in\mathcal{G}_{k-1}$, then $\Upsilon_K(t)$ has no singularities on $(0,\frac{1}{k-1})\cap\mathbb{Q}$.
Thus $\{K\mapsto\frac{1}{2n-1}\Delta\Upsilon'_K(\frac{2}{2n-1})\}_{n=k}^\infty$ gives a homomorphism from $\mathcal{C}/\mathcal{G}_{k-1}$ to $\mathbb{Z}^\infty$.
Hence $\{J_n'\}_{n=k}^\infty$ form a basis for a $\mathbb{Z}^\infty$ summand $\mathcal{C}/\mathcal{G}_{k-1}$.
Note that since $J_n'$'s are topologically slice, the above argument can be restricted to the subgroup
$\mathcal{C}_{TS}/(\mathcal{C}_{TS}\cap\mathcal{G}_{k-1})$ of $\mathcal{C}/\mathcal{G}_{k-1}$ to complete the proof.
\hfill$\Box$

\end{document}